\documentclass{amsart}
\usepackage{amsmath,amsfonts,amssymb,amscd,amsthm,epsf}
\usepackage{pinlabel}
\usepackage{hyperref}
\usepackage{caption}
\newtheorem{prop}{Proposition}[section]
\newtheorem{thm}[prop]{Theorem}

\theoremstyle{definition}

\theoremstyle{remark}
\newtheorem{Remarks}[prop]{Remarks}             
\def\C{{\mathbb C}}

\def\Z{{\mathbb Z}}
\def\R{{\mathbb R}}
\def\Q{{\mathbb Q}}

\def\co{\colon\thinspace}

\begin{document}
\title{On sections of maps from 4-manifolds to the 2-sphere}
\author{Robert E. Gompf}
\address{The University of Texas at Austin}
\email{gompf@math.utexas.edu}
\thanks{The author thanks \.Inan\c c Baykur for helpful comments on the first draft.}

\begin{abstract} 
This note exhibits singular fibrations over the 2-sphere whose regular fibers are connected surfaces of arbitrarily high genus, but which admit no sections. These include achiral Lefschetz fibrations, as well as generic maps for which some disks cannot lift to sections -- in fact, neither hemisphere lifts.
\end{abstract}
\maketitle


\section{Introduction}\label{Intro}

In recent decades, much work has been done to understand various types of smooth maps $f$ from a closed, oriented 4-manifold $X$ to the 2-sphere (for example \cite{M}, \cite{ADK}, \cite{B}, \cite{W}, \cite{GK}, \cite{BS}, and \cite{S} in arbitrary dimensions). These are often treated as singular fibrations of various sorts. Examples include nonconstant holomorphic maps and generic surjections. Of particular interest is the question of when such maps admit sections, i.e., maps from the sphere back to $X$ whose composite with $f$ is the identity on the sphere. For example, when $f$ is a fiber bundle projection, it admits a section unless the fibers have genus 1, and regardless of genus, sections exist over any disk in the sphere. In contrast, this note exhibits examples whose critical locus in $X$ is finite, but which do not admit sections, and generic maps for which the sphere is covered by disks that cannot lift to sections. To avoid the genus-1 exception, we arrange all regular fibers to have arbitrarily high genus. We also require the fibers to be nonempty and connected, to avoid other easy counterexamples. (For such a counterexample, start with a map over the equator that factors through a covering of the circle, extend over both hemispheres by coning, and perturb to be smooth and generic.) To put these results in context, we now consider how the behavior of singular fibrations varies with the dimension of the critical locus.

When the critical locus is empty, we have a fiber bundle. There are two such bundles whose fibers have genus 0, and these both have sections. A simple example whose fibers have genus 1 is the composite $S^3\times S^1\to S^3\to S^2$ of projection followed by the Hopf fibration. This has no section since the Hopf fibration does not. To see the failure explicitly, note that the bundle is trivial over each hemisphere, uniquely up to isotopy. The unique section over the northern hemisphere, restricted to the equator, can be uniquely interpreted in the trivialization over the southern hemisphere as a loop in the torus. Since this loop is homotopically nontrivial, it cannot be extended over the southern hemisphere. In contrast, when the fiber genus is at least 2, the two corresponding trivializations fit together to trivialize the entire bundle, since the diffeomorphism group of the fiber has trivial fundamental group. Thus, there must be a section.

When the critical locus in $X$ is finite, there is a section over any disk. To construct such a section, first note that each fiber $f^{-1}(p)$ has regular points. (Otherwise, the end of $X-f^{-1}(p)$, being simply connected, would lift as a fiber bundle over the $\Z$-cover of the end of $S^2-\{p\}$, contradicting compactness of the fibers.) Choose a regular point of $f$ in each singular fiber (preimage of a critical value). By the Local Submersion Theorem, this gives a local section near each critical value. Since we assume the regular fibers are connected, these local sections can be joined to get a section over any preassigned disk. Thus, the only obstruction to a global section over the sphere is whether we can extend over a final small disk of regular values. This is a difficult problem in general, partly because the section over the large disk need not be unique. In the case of {\em Lefschetz fibrations}, whose critical points have complex quadratic local models, it is still unknown whether sections always exist (outside the previous case of nonsingular torus bundles). However, for {\em achiral Lefschetz fibrations}, where the critical points are complex quadratic in local charts on $X$ that may be orientation-reversing, sections need not exist. Such an example with fibers of genus 1 on $X=S^4$ was constructed by Matsumoto \cite{M}. This admits no section since the intersection pairing vanishes. For higher genus examples, the following observation was informally disseminated by the author several decades ago, but may deserve renewed attention:

\begin{thm}\label{achiral}
There are achiral Lefschetz fibrations over the sphere, with connected fibers of arbitrarily high genus, that do not admit sections.
\end{thm}

For a generic map $f$, the critical locus $K$ in $X$ is a smooth 1-manifold. (In contrast, entire 2-dimensional fibers can lie in $K$ in the holomorphic setting. This happens in elliptic surfaces, for example, and can obstruct sections over disks.) Generic maps have been intensively studied in the past two decades (eg.\ \cite{GK}, \cite{BS}, \cite{S}). Every smooth map $f\co X\to S^2$ is homotopic to such a map for which the fibers are connected and $f(K)$ is an embedded cusped circle \cite[Corollary~1]{W} (see also \cite[Corollary~6.2]{BS}). These frequently have high fiber genus and no sections. For example, there cannot be a section if $f^*[S^2]$ vanishes in $H^2(X;\Q)$. Like the case of finite $K$, however, fiber-connected maps with such $f(K)$ admit sections over the complement of any regular value. This raises the more subtle question of whether sections of fiber-connected generic maps always exist over preassigned disks. Such sections can be useful when they exist. However, at least one submitted manuscript has collapsed over an unsupported existence assumption. We show that sections over disks need not exist by supplying simple counterexamples.

\begin{thm}\label{disk}
There is a generic map $f\co X\to S^2$, whose regular fibers are connected surfaces of arbitrarily high genus, for which neither hemisphere can be lifted to a section.
\end{thm}

\noindent The map does have a section over the complement of both poles, but such sections cannot be extended over either pole.

The critical locus $K$ of a generic map can have several types of points, requiring various local models. However, our examples are built with only one simple type, sometimes called an {\em indefinite fold}. Their local model is
$$f(t,x,y,z)=(t,x^2+y^2-z^2).$$
Thus, the first coordinate is preserved, and in the normal directions we see a Morse critical point of index 1 (or 2 if we rotate by flipping the signs of $t$, $x$ and both image coordinates). Then $K$ is a 1-manifold in $X$ with $f|K$ a smooth immersion. Thus, we avoid the cusps mentioned previously (points where $K$ is tangent to the fibers). When a generic path in $S^2$ crosses $f(K)$, the fibers over it change by surgery on an embedded $S^0$ or $S^1$, depending on the direction of the path. Thus the fiber genus is variable, but in our examples it can be chosen to take only the values $g+1,g,g-1$ for any preassigned $g\ge1$. Furthermore, our critical image $f(K)$ is transverse to the equator, whose preimage is then a 3-manifold with a circle-valued Morse function.

The proofs of these two theorems are independent and comprise the two remaining sections.

\section{Proof of Theorem~\ref{achiral}}\label{Achiral}

To construct an achiral Lefschetz fibration without sections, begin with a generic pencil of curves of a fixed degree $d\ge 2$ on $\C P^2$. This is a {\em Lefschetz pencil}. That is, it has the structure of a Lefschetz fibration everywhere except at the $d^2$ points of the {\em base locus}, where all fibers intersect as near the origin in the local model $\C^2-\{0\}\to\C P^1=S^2$ given by projectivization $(z,w)\mapsto z/w$ (so the fibers are given by complex lines through the origin). Each fiber represents $d$ times a generator of $H_2(\C P^2)\cong\Z$. Let $X_0$ be obtained from $\C P^2$ by deleting a small open ball around each point of the base locus. Then $X_0$ inherits a Lefschetz fibration $f_0\co X_0\to S^2$ whose restriction to each of the  $d^2$ boundary components is the Hopf fibration on $S^3$. Let $X$ be the double of $X_0$, obtained by rounding corners of $\partial (I\times X_0)$, and let $f\co X\to S^2$ be the composite of $f_0$ with the obvious projection. Then $X$ is diffeomorphic to the connected sum $\C P^2\#\overline{\C P^2}\#(d^2-1)S^3\times S^1$. Because the second summand has reversed orientation, each critical point of $f_0$ corresponds to a pair of oppositely oriented critical points of $f$, so the latter is an achiral Lefschetz fibration. Each fiber of $f$ is the double of the corresponding fiber of $f_0$, so it represents $d$ times the sum of generators of the first two summands. In particular, its intersection number with any other class must be divisible by $d\ge2$. Since any section would have intersection number 1 with the fibers, sections cannot exist.

\begin{Remarks}
a) It is easily checked that the regular fibers of $f$ have genus $(d-1)(d-2)+d^2-1=(d-1)(2d-1)$ and that the critical locus has $3(d-1)^2$ quadratic singularities of each orientation.
\item[b)] Another way to describe this construction is to first blow up the $d^2$ points of the base locus of the pencil, obtaining a Lefschetz fibration on $\C P^2\# d^2\overline{\C P^2}$ with $d^2$ disjoint sections of square $-1$. Then $X$ is obtained from two oppositely oriented copies of this by normal connected sum along the corresponding pairs of sections. From this perspective, we can obtain achiral Lefschetz fibrations much more generally by summing oppositely oriented Lefschetz fibrations along pairs of sections that have the same square before reversing orientation. The resulting 4-manifolds will tend to have simple diffeomorphism types (cf.\ \cite{G}); for example, if both summands have spheres of negative square avoiding the sections then the Seiberg-Witten invariants of the result will vanish with both orientations. However, one might at least expect to realize other fiber genera for achiral Lefschetz fibrations without sections by starting with other Lefschetz pencils with homologically nonprimitive fibers.
\end{Remarks}

\section{Proof of Theorem~\ref{disk}}\label{Disk}

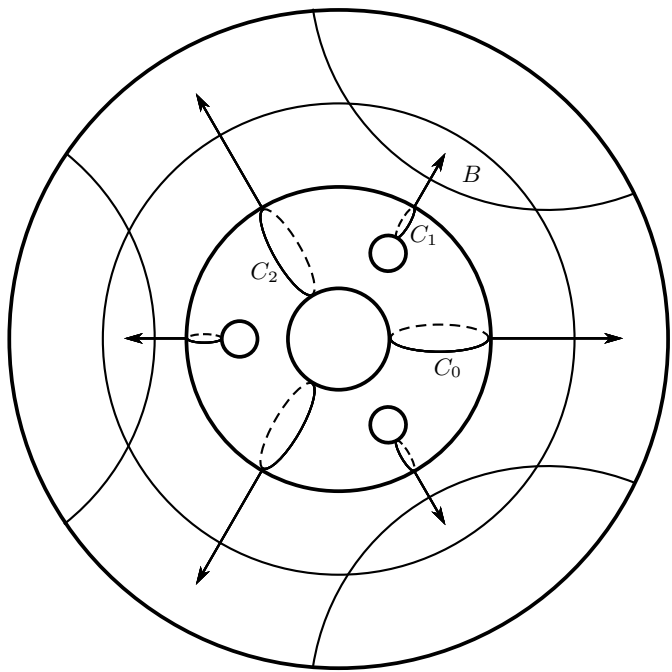
\begin{figure}
\labellist
\small\hair 2pt
\pinlabel {$C_0$} at 166 114
\pinlabel {$C_1$} at 157 164
\pinlabel {$C_2$} at 97 150
\pinlabel {$B$} at 175 187
\endlabellist
\centering
\includegraphics{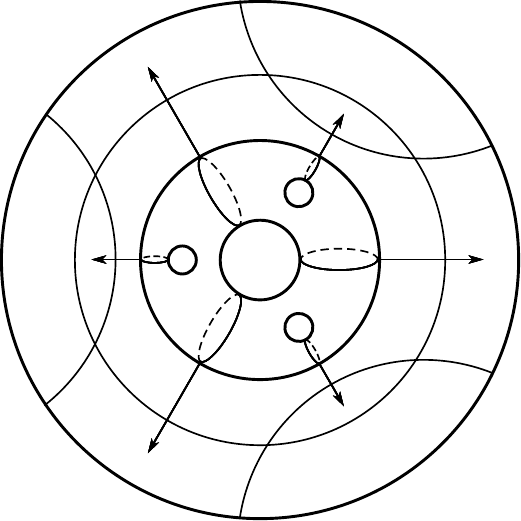}
\caption{Schematic of a generic map to a disk admitting no section.}
\label{fig}
\end{figure}

We begin by constructing a generic map $f_0\co (X_0,\partial X_0)\to (D,\partial D)$ from a compact 4-manifold to the unit disk centered at the origin in $\R^2$, such that $f_0$ has no section, and has boundary-transverse critical locus consisting of indefinite folds. For a preassigned $g\ge1$, each regular fiber will have genus $g+1, g$ or $g-1$. Figure~\ref{fig} represents the construction by the usual schematic when $g=3$: The surface in the center represents the fiber $F$ over the center point of $D$, with genus $g+1$ and an obvious $\Z/g$-action. The background shows $D$ (bounded by the outer circle) with critical image $f_0(K)$ given by the fine circle and $g$ arcs embedded rel boundary. The arrows indicate which circles of $F$ are surgered out as we cross $f_0(K)$ as shown. For example, when we traverse the ray $\theta=0$ (in polar coordinates) from the origin to $\partial D$, crossing $f_0(K)$ surgers out the curve $C_0$ to create fibers of genus $g$. For rays obtained by gradually increasing $\theta$, we eventually encounter a second critical value corresponding to surgering $C_1$. (Since $C_0$ and $C_1$ are disjoint, we can do the surgeries consecutively as indicated, to get fibers of genus $g-1$ near that part of $\partial D$.) As $\theta$ increases, the ray passes the first intersection point of $f_0(K)$, reversing the order of the two surgeries. On the bigon $B$ between these two surgeries, we can interpret $C_0$ as a circle in the genus-$g$ surface obtained by surgering only $C_1$. After this first surgery, $C_0$ is isotopic to the circle $C_2$, so we equivalently have consecutive surgeries on $C_1$ and $C_2$. As $\theta$ increases further, the ray passes another intersection, again reversing the order of these surgeries. Then the surgery on $C_1$ disappears at the end of that arc, and we are left with just a surgery on $C_2$. The net effect is that we have moved $C_0$ to $C_2$ by temporarily surgering out the obstructing tube. Continue in this way around $F$ until $\theta$ reaches $2\pi$. To turn the schematic into an actual 4-manifold $X_0$ and map $f_0$, start with projection of $F\times D^2$ to a small central disk in $D$. As the radius of the disk increases, we encounter the arcs of critical values. Each is realized by attaching a 0-framed 2-handle to the corresponding circle in the appropriate fiber. Then we encounter the critical circle, which attaches a round 2-handle $S^1\times D^2\times I$ to the new boundary.

To see that this map $f_0$ has no section, suppose that it does have a section $\sigma$, that is, a smoothly varying choice of a point in each fiber. Consider its restriction to $\partial D$. The nonsingular fibers along $\partial D$ all have genus $g$ or $g-1$. At $\theta=0$ the fiber is obtained from $F$ by surgering out $C_0$. As $\theta$ increases to $2\pi$, this surgered circle sweeps all the way around $F$, with obstructing tubes temporarily vanishing as needed. This must sweep $\sigma$ ahead of it, so that $\sigma|\partial D$ also travels around the fiber. When we consider $\sigma$ over circles of decreasing radius, the surgeries disappear and we are left with a loop in $F$. Since this loop is homotopically essential, $\sigma$ cannot be extended over all of $D$ (although we have essentially constructed a section over $D-\{0\}$).

To finish the proof, simply double $D$ and $X_0$ to get a fibration $f\co X\to S^2$. This has a section over the complement of the two poles, but no section over either hemisphere.

\begin{Remarks}
a) The manifold $X$ is $\# (g+1)S^3\times S^1$. To see this, decompose $F\times D^2$ as a handlebody with a unique 2-handle $h$, whose cocore $\delta$ is $\{p\}\times D^2$. The above description now exhibits $X_0$  as a handlebody once we decompose the round 2-handle as a 2-handle and 3-handle. The attaching sphere of the 3-handle intersects $\delta$ transversely in a unique point (where $C_0$ sweeps over $p$ as it traverses the model fiber), so the 3-handle cancels $h$. The other 2-handles of $X_0$ cancel 1-handles of $F\times D^2$, so $X_0$ is a 0-handle with $g+1$ 1-handles attached, and its double $X=\partial(I\times X_0)$ is as claimed. This shows again that $f_0$ has no section, which would give a class in $H_2(X_0,\partial X_0)=0$ pairing nontrivially with a fiber.

\item[b)] For a variant with $H_2\ne 0$, replace the bundle over $D$ by its pullback via the branched covering map $z\mapsto z^d$. This preserves the model fiber, but now the figure has $gd$ disjoint critical arcs along with the critical circle. This results in $d$ 2-handles attached along each $C_i$ with $i$ odd, contributing $g(d-1)$ generators to $H_2(X_0)$ and $g(d-1)$ hyperbolic pairs to the intersection pairing of $X$. The 3-handle no longer cancels, but the previous analysis shows the fiber has order $d$. Either previous method still shows there is no section over $D$, although the resulting $f$ still has a section over the complement of both poles.
\end{Remarks}


\begin{thebibliography}{MM}

\bibitem[ADK05]{ADK}
D.\ Auroux, S.\ Donaldson and L.\ Katzarkov, 
{\em Singular Lefschetz pencils}, 
Geom.\ Topol.\ {\bf 9} (2005), 1043--1114.

\bibitem[B08]{B}
R.\ \.I.\ Baykur, 
{\em Existence of Broken Lefschetz Fibrations}, 
Int.\ Math.\ Res.\ Not.\ {\bf 2008} (2008), Article ID rnn101, 15 pages, doi:10.1093/imrn/rnn101

\bibitem[BS23]{BS}
R.\ \.I.\ Baykur and O.\  Saeki, 
{\em Simplifying indefinite fibrations on 4-manifolds}, 
Trans.\ Amer.\ Math.\ Soc.\ {\bf 376} (2023), 3011--3062.

\bibitem[GK15]{GK}
D.\ Gay and R.\ Kirby, 
{\em Indefinite Morse 2-functions; broken fibrations and generalizations}, 
Geom.\ Topol.\ {\bf 19} (2015), 2465--2534.

\bibitem[G88]{G}
R.\ Gompf,
{\em On sums of algebraic surfaces}, 
Invent.\ Math.\ {\bf 94} (1988), 171--174.

\bibitem[M82]{M}
Y.\ Matsumoto, 
{\em On 4-manifolds fibered by tori}, 
Proc.\ Japan Acad.\ {\bf 58} (1982), 298--301.

\bibitem[S25]{S}
O.\ Saeki, 
{\em Simplifying generic smooth maps to the 2-sphere and to the plane}, 
Preprint (revised 2025), arXiv:2407.10145

\bibitem[W10]{W}
J.\ Williams, 
{\em The h-principle for broken Lefschetz fibrations}, 
Geom.\ Topol.\ {\bf 14}, no.\ 2 (2010), 1015--1061.

\end{thebibliography}
\end{document}